\documentclass{birkmult}

\usepackage[T1]{fontenc}

\usepackage[pdftex,breaklinks=true,colorlinks=true,pdfstartview=FitH,       urlcolor=blue]{hyperref}

\theoremstyle{plain}
\newtheorem{theorem}{Theorem}
\newtheorem{lemma}{Lemma}

\theoremstyle{remark}
\newtheorem{example}{Example}

\numberwithin{equation}{section}

\DeclareMathOperator{\RE}{Re}
\DeclareMathOperator{\IM}{Im}

\begin{document}
%-------------------------------------------------------------------------
% editorial commands: to be inserted by the editorial office
%
%\firstpage{1}
%\volume{228}
%\DOI{003-0001}
%\seriesextraline{This is the Concrete Title of this Book\br H.E. R and S.T.C. W, Eds.}
%\firstpage{1}
\title[Series representations of the multiparameter fBm] {A family of series representations of the multiparameter fractional Brownian motion}

\author{Anatoliy Malyarenko}

\address{%
Division of Applied Mathematics\\
School of Education, Culture and Communication\\
Mälardalen University\\
SE 721 23 Västerås, Sweden}

\email{anatoliy.malyarenko@mdh.se}

\thanks{This work is supported by the Swedish Institute grant SI--01424/2007.}

\subjclass{Primary 60G60; Secondary 33C60}

\keywords{Multiparameter fractional Brownian motion, series representation, Meijer $G$-function.}

\date{\today}

\begin{abstract}
We derive a family of series representations of the multiparameter fractional Brownian motion in the centred ball of radius~$R$ in the $N$-dimensional space $\mathbb{R}^N$. Some known examples of series representations are shown to be the members of the family under consideration.
\end{abstract}

%%% ----------------------------------------------------------------------
\maketitle
%%% ----------------------------------------------------------------------
%\tableofcontents

\section{Introduction}

The \emph{fractional Brownian motion} with Hurst parameter
$H\in(0,1)$ is defined as the centred Gaussian process $\xi(t)$ with the autocorrelation function
\[
R(s,t)=\mathsf{E}\xi(s)\xi(t)=\frac{1}{2}(|s|^{2H}+|t|^{2H}-|s-t|^{2H}).
\]
This process was introduced by Kolmogorov (1940) and became a popular statistical model after the paper by Mandelbrot and van Ness (1968).

There exist two multiparameter extensions of the fractional Brownian motion. Both extensions are centred Gaussian random fields on the space $\mathbb{R}^N$. The \emph{multiparameter fractional Brownian sheet} has the autocorrelation function
\[
R(\mathbf{x},\mathbf{y})=\frac{1}{2^N}\prod^N_{j=1}(|x_j|^{2H_j}+|y_j|^{2H_j}-|x_j-y_j|^{2H_j}),\qquad
H_j\in(0,1),
\]
while the \emph{multiparameter fractional Brownian motion} has the autocorrelation function
\begin{equation}\label{auto}
R(\mathbf{x},\mathbf{y})=\frac{1}{2}(\|\mathbf{x}\|^{2H}+\|\mathbf{y}\|^{2H}
-\|\mathbf{x}-\mathbf{y}\|^{2H}),
\end{equation}
where $\|\cdot\|$ denote the Euclidean norm in $\mathbb{R}^N$ and where $H\in(0,1)$.

Malyarenko (2008) derived a series expansion of the multiparameter fractional Brownian motion. His expansion converges almost surely (a. s.) in the space $C(\mathcal{B})$ of continuous functions in the centred ball $\mathcal{B}=\{\,\mathbf{x}\in\mathbb{R}^N\colon\|\mathbf{x}\|\leq 1\,\}$.

In fact, the above mentioned series expansion is a member of a family of series expansions.
To describe this family, introduce the following notation. Let
\begin{equation}\label{constant}
c_{NH}=\sqrt{2\pi^{(N-2)/2}\Gamma(N/2+H)\Gamma(H+1)\sin(\pi H)},
\end{equation}
where $\Gamma$ denote the gamma function. A review of special functions
is given in Section~\ref{special}. Denote
\begin{equation}\label{kernels}
a_m(s,u)=
\begin{cases}
c_{NH}u^{H-1/2}G^{2,0}_{2,2}\left(\dfrac{u^2}{s^2}\left|
\begin{array}{l}
N/2,1\\
0,1-H
\end{array}
\right.\right),&m=0,\\
c_{NH}s^{2H-m}u^{m-H-1/2}G^{1,0}_{1,1}\left(\dfrac{u^2}{s^2}\left|
\begin{array}{l}
N/2+H\\
0
\end{array}
\right.\right),&m\geq 1,
\end{cases}
\end{equation}
where $G^{m,n}_{p,q}\left(z\left|
\begin{array}{l}
a_1,\dots,a_n,a_{n+1},\dots,a_p\\
b_1,\dots,b_m,b_{m+1},\dots,b_q
\end{array}
\right.\right)$ is the Meijer $G$-function.

Let $\mathbb{Z}_+$ be the set of all nonnegative integers. For a fixed $m\in\mathbb{Z}_+$, there exist
\[
h(m,N)=\frac{(2m+N-2)(m+N-3)!}{(N-2)!m!}
\]
different real-valued spherical harmonics $S^l_m(\mathbf{x}/\|\mathbf{x}\|)$ of degree~$m$.
Fix $R>0$, and for each pair $(m,l)$ with $m\in\mathbb{Z}_+$ and $1\leq l\leq h(m,N)$,
let $\{\,e^l_{mn}(u)\colon n\geq 1\,\}$ be a basis in the Hilbert space $L^2[0,R]$.
Let $b^l_{mn}(s)$ be the Fourier coefficients of the function $a_m(s,u)$
with respect to the introduced basis:
\begin{equation}\label{coefficients}
b^l_{mn}(s)=\int^R_0a_m(s,u)e^l_{mn}(u)\,du.
\end{equation}
Finally, let $\{\,\xi^l_{mn}\colon m\in\mathbb{Z}_+, n\geq 1, 1\leq l\leq h(m,N)\,\}$ be the set of independent standard normal random variables.

\begin{theorem}\label{th1}
For any choice of the bases $\{\,e^l_{mn}(u)\colon n\geq 1\,\}$, the
multiparameter fractional Brownian motion $\xi(\mathbf{x})$ has the
following series expansion
\begin{equation}\label{expansion}
\xi(\mathbf{x})=\sum^{\infty}_{m=0}\sum^{h(m,N)}_{l=1}\sum^{\infty}_{n=1}
b^l_{mn}(\|\mathbf{x}\|)S^l_m(\mathbf{x}/\|\mathbf{x}\|)\xi^l_{mn}.
\end{equation}
The series \eqref{expansion} converges in mean square in the centred closed ball $\mathcal{B}_R=\{\,\mathbf{x}\in\mathbb{R}^N\colon\|\mathbf{x}\|\leq R\,\}$.
\end{theorem}

In Section~\ref{special} we review the necessary definitions and properties of some special functions. In Section~\ref{proof} we give an outline of proof of Theorem~\ref{th1}. Section~\ref{examples} contains examples, while proofs of technical lemmas are postponed to Section~\ref{technical}.

\section{Special functions}\label{special}

This section is intended for readers who are not experts in the theory of special functions. We review definitions and state some elementary properties of the relevant special functions.

This material can be found online at \href{http://functions.wolfram.com}{\texttt{http://functions.wolfram.com}}, a comprehensive online compendium of formulas.

\subsection{The gamma function}

The \emph{gamma function} of a complex variable $z$ with $\RE z>0$ is defined by\footnote{\href{http://functions.wolfram.com/06.05.02.0001.01}{\texttt{http://functions.wolfram.com/06.05.02.0001.01}}}
\[
\Gamma(z)=\int^{\infty}_0t^{z-1}e^{-t}\,dt.
\]
By partial integration, we obtain\footnote{\href{http://functions.wolfram.com/06.05.16.0004.01}{\texttt{http://functions.wolfram.com/06.05.16.0004.01}}} %
\begin{equation}\label{recurrent}
\Gamma(z-1)=\frac{\Gamma(z)}{z-1}.
\end{equation}
This formula is used to extend $\Gamma$ to an analytic function of  $z\in\mathbb{C}\setminus\mathbb{Z}_-$, where $\mathbb{Z}_-=\{0,-1,\dots,-n,\dots\}$. The points $z\in\mathbb{Z}_-$ are the poles.

\subsection{The Meijer $G$-function}

Let $m$, $n$, $p$, and $q$ be four nonnegative integers with $0\leq m\leq q$ and $0\leq n\leq p$. Let $a_1$, \dots, $a_p$, $b_1$, \dots, $b_q$ be points in the complex plane. Assume that for each $k=1$, $2$, \dots, $p$ and for each $i=1$, $2$, \dots, $q$ we have $a_k-b_i+1\notin\mathbb{Z}_+$. Then, there exists an infinite contour $\mathcal{L}$ that separates the poles of $\Gamma(1-a_k-s)$ at $s=1-a_k+j$, $j\in\mathbb{Z}_+$ from the poles of $\Gamma(b_i+s)$ at $s=-b_i-l$, $l\in\mathbb{Z}_+$. The \emph{Meijer $G$-function} is defined by\footnote{\href{http://functions.wolfram.com/07.34.02.0001.01}{\texttt{http://functions.wolfram.com/07.34.02.0001.01}}} \[
\begin{aligned}
&G^{m,n}_{p,q}\left(z\left|
\begin{array}{l}
a_1,\dots,a_n,a_{n+1},\dots,a_p\\
b_1,\dots,b_m,b_{m+1},\dots,b_q
\end{array}
\right.\right)=\\
&\quad\frac{1}{2\pi i}\int_{\mathcal{L}}\frac{\Gamma(s+b_1)\cdots\Gamma(s+b_m)\Gamma(1-a_1-s)\cdots\Gamma(1-a_n-s)}
{\Gamma(s+a_{n+1})\cdots\Gamma(s+a_p)\Gamma(1-b_{m+1}-s)\cdots\Gamma(1-b_q-s)}z^{-s}\,ds.
\end{aligned}
\]
The number $p+q$ is the \emph{order} of the Meijer $G$-function.

The classical Meijer's integral from two $G$-functions\footnote{\href{http://functions.wolfram.com/07.34.21.0011.01}{\texttt{http://functions.wolfram.com/07.34.21.0011.01}}} is:
\begin{equation}\label{classic}
\begin{aligned}
&\int^{\infty}_0\tau^{\alpha-1}G^{s,t}_{u,v}\left(w\tau\left|
\begin{array}{l}
c_1,\dots,c_t,c_{t+1},\dots,c_u\\
d_1,\dots,d_s,d_{s+1},\dots,d_v
\end{array}
\right.\right)\\
&\quad\times G^{m,n}_{p,q}\left(z\tau\left|
\begin{array}{l}
a_1,\dots,a_n,a_{n+1},\dots,a_p\\
b_1,\dots,b_m,b_{m+1},\dots,b_q
\end{array}
\right.\right)\,d\tau\\
&\quad=w^{-\alpha}G^{m+t,n+s}_{v+p,u+q}\left(\frac{z}{w}\left|
\begin{array}{l}
a_1,\dots,a_n,1-\alpha-d_1,\dots,1-\alpha-d_v,a_{n+1},\dots,a_p\\
b_1,\dots,b_m,1-\alpha-c_1,\dots,1-\alpha-c_u,b_{m+1},\dots,b_q
\end{array}
\right.\right).
\end{aligned}
\end{equation}

A plenty of special functions are specialised values of the Meijer $G$-function.
In particular, the \emph{Gegenbauer polynomials} $C^{\lambda}_n(x)$ appear
as\footnote{\href{http://functions.wolfram.com/07.34.03.0105.01}{\texttt{http://functions.wolfram.com/07.34.03.0105.01}}} %
\begin{equation}\label{Gegenbauer}
\begin{aligned}
G^{0,2}_{2,2}\left(z\left|
\begin{array}{l}
a,c\\
b,b+1/2
\end{array}
\right.\right)&=\frac{\Gamma(2b-2c+2)\vartheta(|z|-1)}{\Gamma(a-2b+c-1/2)(2(a-2b+c-1))_{2b-2c+1}}\\
&\quad\times z^b(z-1)^{a-2b+c-3/2}C^{a-2b+c-1}_{2b-2c+1}(\sqrt{z}),
\end{aligned}
\end{equation}
where $(a)_n=\Gamma(a+n)/\Gamma(a)$ is the Pochhammer symbol, and
\[
\theta(x)=
\begin{cases}
1,&x\geq 0,\\
0,&x<0
\end{cases}
\]
is the unit step function. In particular, the \emph{Legendre polynomials} $P_n(x)$ are just Gegenbauer polynomials with $\lambda=0$.

The \emph{generalised hypergeometric function} ${}_pF_q(a_1,\dots,a_p;b_1,\dots,b_q;z)$ appears as\footnote{\href{http://functions.wolfram.com/07.31.26.0004.01}{\texttt{http://functions.wolfram.com/07.31.26.0004.01}}} \[
{}_qF_p(a_1,\dots,a_p;b_1,\dots,b_q;z)=\frac{\Gamma(b_1)\cdots\Gamma(b_q)}{\Gamma(a_1)\cdots\Gamma(a_p)}
G^{1,p}_{p,q+1}\left(-z\left|
\begin{array}{l}
1-a_1,\dots,1-a_p\\
0,1-b_1,\dots,1-b_q
\end{array}
\right.\right).
\]

The \emph{Bessel functions} appear as\footnote{\href{http://functions.wolfram.com/03.01.26.0107.01}{\texttt{http://functions.wolfram.com/03.01.26.0107.01}}} %
\begin{equation}\label{bessel}
J_{\nu}(z)=G^{1,0}_{0,2}\left(\frac{z^2}{4}\left|
\begin{array}{l}
\boldsymbol{\cdot} \\
\nu/2,-\nu/2
\end{array}
\right.\right).
\end{equation}

\subsection{Spherical harmonics}

Let $m$ be a nonnegative integer, and let $m_0$, $m_1$, \dots, $m_{N-2}$ be integers satisfying the following condition
\[
m=m_0\geq m_1\geq\dots\geq m_{N-2}\geq 0.
\]
Let $\mathbf{x}=(x_1, x_2,\dots, x_N)$ be a point in the space $\mathbb{R}^N$. Let
\[
r_k=\sqrt{x^2_{k+1}+x^2_{k+2}+\dots+x^2_N},
\]
where $k=0$, $1$, \dots, $N-2$. Consider the following functions
\[
\begin{aligned}
H(m_k,\pm,\mathbf{x})&=\left(\frac{x_{N-1}+ix_N}{r_{N-2}}\right)^{\pm m_{N-2}}r^{m_{N-2}}_{N-2}\prod^{N-3}_{k=0}r^{m_k-m_{k+1}}_k\\
&\quad\times C^{m_{k+1}+(N-k-2)/2}_{m_k-m_{k+1}}\left(\frac{x_{k+1}}{r_k}\right),\\
\end{aligned}
\]
and denote
\[
Y(m_k,\pm,\mathbf{x})=r^{-m}_0H(m_k,\pm,\mathbf{x}).
\]
The functions $Y(m_k,\pm,\mathbf{x})$ are called the (complex-valued) \emph{spherical harmonics}. They are orthogonal in the Hilbert space $L^2(S^{N-1})$ of the square integrable functions on the unit sphere $S^{N-1}$, and the square of the length of the vector $Y(m_k,\pm,\mathbf{x})$ is
\[
L(m_k)=2\pi\prod^{N-2}_{k=1}\frac{\pi2^{k-2m_k-N+2}\Gamma(m_{k-1}+m_k+N-1-k)}
{(m_{k-1}+(N-1-k)/2)(m_{k-1}-m_k)![\Gamma(m_k+(N-1-k)/2)]^2}.
\]

Let $l=l(m_k,\pm)$ be the number of the symbol $(m_0, m_1, \dots, m_{N-2},\pm)$ in the lexicographic ordering. The \emph{real-valued spherical harmonics}, $S^l_m(\mathbf{x})$, can be defined as
\[
S^l_m(\mathbf{x})=
\begin{cases}
Y(m_k,+,\mathbf{x})/\sqrt{L(m_k)},&m_{N-2}=0,\\
\sqrt{2}\RE Y(m_k,+,\mathbf{x})/\sqrt{L(m_k)},&m_{N-2}>0,l=l(m_k,+),\\
-\sqrt{2}\IM Y(m_k,-,\mathbf{x})/\sqrt{L(m_k)},&m_{N-2}>0,l=l(m_k,-).
\end{cases}
\]
The spherical harmonics $S^l_m(\mathbf{x})$ form a basis in the Hilbert space $L^2(S^{N-1})$.

\section{Proof of Theorem~1 modulo technical lemmas}\label{proof}

Note that the multiparameter fractional Brownian motion $\xi(\mathbf{x})$ is \emph{weakly isotropic}, i.e., the autocorrelation function \eqref{auto} is invariant with respect to the group $O(N)$ of the orthogonal matrices of order~$N$. Let $t>0$, let $S_t=\{\,\mathbf{x}\in\mathbb{R}^N\colon\|\mathbf{x}\|=t\,\}$ be the centred sphere in the space $\mathbb{R}^N$, and let $d\omega$ be the Lebesgue surface measure on $S_t$. Let $\eta(\mathbf{x})$ be a centred weakly isotropic random field. In fact, the autocorrelation function $R(\mathbf{x},\mathbf{y})$ of the random field $\eta(\mathbf{x})$ is a function $R(s,t,u)$ of the three real variables $s=\|\mathbf{x}\|$, $t=\|\mathbf{y}\|$, and $u$ being the cosine of the angle between the vectors $\mathbf{x}$ and $\mathbf{y}$. Yadrenko (1983) proved that the stochastic processes
\[
X^l_m(t)=\int_{S_t}\eta(\mathbf{x})S^l_m(\mathbf{x}/\|\mathbf{x}\|)\,d\omega
\]
are centred and uncorrelated. The autocorrelation function of the process $X^l_m(t)$ is
\[
R_m(s,t)=c\int^1_{-1}R(s,t,u)C^{(N-2)/2}_m(u)(1-u^2)^{(N-3)/2}\,du,
\]
with
\[
c=\frac{2^{N-2}\pi^{(N-2)/2}m!\Gamma((N-2)/2)}{\Gamma(m+N-2)}.
\]
The random field $\eta(\mathbf{x})$ can be represented as
\begin{equation}\label{presentation}
\eta(\mathbf{x})=\sum^{\infty}_{m=0}\sum^{h(m,N)}_{l=1}X^l_m(\|\mathbf{x}\|)
S^l_m(\mathbf{x}/\|\mathbf{x}\|),
\end{equation}
where the series converges in mean square for any $\mathbf{x}\in\mathbb{R}^N$.

The multiparameter fractional Brownian motion $\xi(\mathbf{x})$ is Gaussian random field. Therefore, the corresponding stochastic processes $X^l_m(t)$ are Gaussian and independent.

\begin{lemma}\label{l1}
The autocorrelation function of the stochastic process $X^l_m(t)$ has the form
\begin{equation}\label{xlmauto}
R_m(s,t)=
\begin{cases}
\frac{c_{NH}^2}{2}s^{2H}G^{2,2}_{4,4}\left(\dfrac{s^2}{t^2}\left|
\begin{array}{l}
1-H,0,N/2,1\\
0,1-H,1-H-N/2,-H
\end{array}
\right.\right),&m=0,\\
\frac{c^2_{NH}}{2}s^mt^{2H-m}G^{1,1}_{2,2}\left(\dfrac{s^2}{t^2}\left|
\begin{array}{l}
H+1-m,N/2+H\\
0,1-N/2-m
\end{array}
\right.\right),&m\geq 1.
\end{cases}
\end{equation}
\end{lemma}

Recall that a function $a(s,u)\colon(0,\infty)\times(0,\infty)\mapsto\mathbb{R}$ is called the \emph{Volterra kernel}, if it is locally square integrable, and
\begin{equation}\label{Volterra}
a(s,u)=0\quad\text{for}\quad s<u.
\end{equation}
A \emph{Volterra process} with Volterra kernel $a(s,u)$ is a centred Gaussian stochastic process $\eta(t)$ with autocorrelation function
\[
R(s,t)=\int^{\min\{s,t\}}_0a(s,u)a(t,u)\,du.
\]

\begin{lemma}\label{l2}
The stochastic processes $X^l_m(t)$ are Volterra processes with Volterra kernels \eqref{kernels}.
\end{lemma}

Volterra processes are important in the theory of stochastic integration with respect to general Gaussian processes, see Decreusefond (2005) and the references herein.

By Lemma~\ref{l2}, the autocorrelation function of the stochastic process $X^l_m(t)$ has the form
\[
R_m(s,t)=\int^{\min\{s,t\}}_0a_m(s,u)a_m(t,u)\,du.
\]
By \eqref{Volterra}, the last display can be rewritten as
\[
R_m(s,t)=\int^R_0a_m(s,u)a_m(t,u)\,du,\qquad s,t\in[0,R].
\]
By definition of a basis in the Hilbert space $L^2[0,R]$, we obtain
\[
R_m(s,t)=\sum^{\infty}_{n=1}b^l_{mn}(s)b^l_{mn}(t),\qquad 1\leq l\leq h(m,N).
\]
It follows that the stochastic process $X^l_m(t)$ itself has the form
\[
X^l_m(t)=\sum^{\infty}_{n=1}b^l_{mn}(t)\xi^l_{mn},
\]
where the series converges in mean square. Substituting this formula to \eqref{presentation}, we obtain \eqref{expansion}.

We conjecture that \eqref{expansion} converges uniformly a.s.

\section{Examples}\label{examples}

\begin{example}\label{ex1}
Let $\nu$ be a real number, and let $j_{\nu,1}<j_{\nu,2}<\dots<j_{\nu,n}<\dots$ be the positive zeros of the Bessel function $J_{\nu}(u)$. For any $\nu>-1$, the \emph{Fourier--Bessel functions}
\[
\varphi_{\nu,n}(u)=\frac{\sqrt{2u}}{J_{\nu+1}(j_{\nu,n})}J_{\nu}(j_{\nu,n}u),\qquad
n\geq 1
\]
form a basis in the space $L^2[0,1]$ (Watson (1995), Section~18.24). By change of variable we conclude that the functions
\[
e^l_{mn}(u)=\frac{\sqrt{2u}}{RJ_{\nu+1}(j_{\nu,n})}J_{\nu}(R^{-1}j_{\nu,n}u),\qquad
n\geq 1
\]
form a basis in the space $L^2[0,R]$.

To calculate $b^l_{mn}(s)$, use \eqref{coefficients}. First, consider the case of $m=0$:
\[
b^l_{0n}(s)=\int^R_0c_{NH}u^{H-1/2}G^{2,0}_{2,2}\left(\dfrac{u^2}{s^2}\left|
\begin{array}{l}
N/2,1\\
0,1-H
\end{array}
\right.\right)\frac{\sqrt{2u}}{RJ_{\nu+1}(j_{\nu,n})}J_{\nu}(R^{-1}j_{\nu,n}u)\,du.
\]
Using \eqref{Volterra}, rewrite this formula as
\[
b^l_{0n}(s)=\frac{c_{NH}\sqrt{2}}{RJ_{\nu+1}(j_{\nu,n})}\int^{\infty}_0u^HJ_{\nu}(R^{-1}j_{\nu,n}u)
G^{2,0}_{2,2}\left(\dfrac{u^2}{s^2}\left|
\begin{array}{l}
N/2,1\\
0,1-H
\end{array}
\right.\right)\,du.
\]
To calculate this integral,  use formula 2.24.4.1 from Prudnikov et al (1990) with particular values of parameters $k=l=1$:
\begin{equation}\label{22441}
\begin{aligned}
&\int^{\infty}_0u^{\alpha-1}J_{\nu}(bu)
G^{m,r}_{p,q}\left(\omega u^2\left|
\begin{array}{l}
a_1,\dots,a_p\\
b_1,\dots, b_q
\end{array}
\right.\right)\,du=\\
&\quad\frac{2^{\alpha-1}}{b^{\alpha}}G^{m,r+1}_{p+2,q}\left(
\frac{4\omega}{b^2}\left|
\begin{array}{l}
1-(\alpha+\nu)/2,a_1,\dots,a_p,1-(\alpha-\nu)/2\\
b_1,\dots,b_q
\end{array}
\right.\right),
\end{aligned}
\end{equation}
where $\alpha=H+1$, $b=R^{-1}j_{\nu,n}$, $m=p=q=2$, $r=0$, $\omega=s^{-2}$, $a_1=N/2$, $a_2=1$, $b_1=0$, and $b_2=1-H$. We obtain
\[
b^l_{0n}(s)=\frac{c_{NH}2^{H+1/2}R^H}{J_{\nu+1}(j_{\nu,n})j_{\nu,n}^{H+1}}G^{2,1}_{4,2}\left(
\frac{4R^2}{s^2j^2_{\nu,n}}\left|
\begin{array}{l}
1-\frac{H+1+\nu}{2},\frac{N}{2},1,1-\frac{H+1-\nu}{2}\\
0,1-H
\end{array}
\right.\right).
\]

Formula 8.2.2.9 from Prudnikov et al (1990) states
\begin{equation}\label{8229}
G^{m,n}_{p,q}\left(z\left|
\begin{array}{l}
a_1,\dots,a_{p-1},b_1\\
b_1,\dots, b_q
\end{array}
\right.\right)=G^{m-1,n}_{p-1,q-1}\left(z\left|
\begin{array}{l}
a_1,\dots,a_{p-1}\\
b_2,\dots, b_q
\end{array}
\right.\right).
\end{equation}
If we put $\nu=1-H$, then \eqref{8229} decreases the order of the Mejer $G$-function from $6$ to $4$. We get
\[
b^l_{0n}(s)=\frac{c_{NH}2^{H+1/2}R^H}{J_{2-H}(j_{1-H,n})j_{1-H,n}^{H+1}}G^{1,1}_{3,1}\left(
\frac{4R^2}{s^2j^2_{1-H,n}}\left|
\begin{array}{l}
0,N/2,1\\
0
\end{array}
\right.\right).
\]
For further simplification, use the symmetry
relation\footnote{\href{http://functions.wolfram.com/07.34.17.0012.01}{\texttt{http://functions.wolfram.com/07.34.17.0012.01}}} %
\begin{equation}\label{symmetry}
G^{m,n}_{p,q}\left(z\left|
\begin{array}{l}
a_1,\dots,a_p\\
b_1,\dots,b_q
\end{array}
\right.\right)=G^{n,m}_{q,p}\left(\frac{1}{z}\left|
\begin{array}{l}
1-b_1,\dots,1-b_q\\
1-a_1,\dots,1-a_p
\end{array}
\right.\right).
\end{equation}
We get
\[
b^l_{0n}(s)=\frac{c_{NH}2^{H+1/2}R^H}{J_{2-H}(j_{1-H,n})j_{1-H,n}^{H+1}}
G^{1,1}_{1,3}\left(
\frac{s^2j^2_{1-H,n}}{4R^2}\left|
\begin{array}{l}
1\\
1,1-N/2,0
\end{array}
\right.\right).
\]
Then, use the argument transformation\footnote{\href{http://functions.wolfram.com/07.34.16.0001.01}{\texttt{http://functions.wolfram.com/07.34.16.0001.01}}} %
\begin{equation}\label{argument}
G^{m,n}_{p,q}\left(z\left|
\begin{array}{l}
\alpha+a_1,\dots,\alpha+a_p\\
\alpha+b_1,\dots,\alpha+b_q
\end{array}
\right.\right)=z^{\alpha}G^{m,n}_{p,q}\left(z\left|
\begin{array}{l}
a_1,\dots,a_p\\
b_1,\dots,b_q
\end{array}
\right.\right)
\end{equation}
with $\alpha=1$ to obtain
\[
b^l_{0n}(s)=\frac{c_{NH}2^{H+1/2}R^H}{J_{2-H}(j_{1-H,n})j_{1-H,n}^{H+1}}
\frac{s^2j^2_{1-H,n}}{4R^2}
G^{1,1}_{1,3}\left(
\frac{s^2j^2_{1-H,n}}{4R^2}\left|
\begin{array}{l}
0\\
0,-N/2,-1
\end{array}
\right.\right),
\]
and use the formula\footnote{\href{http://functions.wolfram.com/07.22.26.0004.01}{\texttt{http://functions.wolfram.com/07.22.26.0004.01}}} \[
{}_1F_2(a_1;b_1,b_2;z)=\frac{\Gamma(b_1)\Gamma(b_2)}{\Gamma(a_1)}G^{1,1}_{1,3}\left(-z\left|
\begin{array}{l}
1-a_1\\
0,1-b_1,1-b_2
\end{array}
\right.\right)
\]
with $a_1=1$, $b_1=N/2+1$, $b_2=2$, and $z=-s^2j^2_{1-H,n}/(4R^2)$ to get
\[
\begin{aligned}
b^l_{0n}(s)&=\frac{c_{NH}2^{H+1/2}R^H}{J_{2-H}(j_{1-H,n})j_{1-H,n}^{H+1}}
\frac{s^2j^2_{1-H,n}}{4R^2}\frac{1}{\Gamma(N/2+1)}\\
&\quad\times{}_1F_2(1;N/2+1,2;-s^2j^2_{1-H,n}/(4R^2)).
\end{aligned}
\]
Finally, use the formula\footnote{\href{http://functions.wolfram.com/07.22.03.0030.01}{\texttt{http://functions.wolfram.com/07.22.03.0030.01}}}
\[
{}_1F_2(1;2,c;z)=I_{c-2}(2\sqrt{z})\Gamma(c)z^{-c/2}+\frac{1-c}{z}
\]
with $c=N/2+1$ and $z=-s^2j^2_{1-H,n}/(4R^2)$. Here $I$ denote the \emph{modified Bessel function}
\[
I_{\nu}(z)=e^{-\nu\pi i/2}J_{\nu}(e^{\pi i/2}z).
\]
We obtain
\begin{equation}\label{case0}
\begin{aligned}
b^l_{0n}(s)&=\frac{c_{NH}2^{H+1/2}R^H}{J_{2-H}(j_{1-H,n})j_{1-H,n}^{H+1}\Gamma(N/2)}\\
&\quad\times\left[2^{(N-2)/2}\Gamma(N/2)\frac{J_{(N-2)/2}(R^{-1}j_{1-H,n}s)}{(R^{-1}j_{1-H,n}s)^{(N-2)/2}}-1\right].
\end{aligned}
\end{equation}

Continue with the remaining case of $m\geq 1$. Using \eqref{coefficients}, we obtain
\[
\begin{aligned}
b^l_{mn}(s)&=\int^R_0c_{NH}s^{2H-m}u^{m-H-1/2}G^{1,0}_{1,1}\left(\dfrac{u^2}{s^2}\left|
\begin{array}{l}
N/2+H\\
0
\end{array}
\right.\right)\\
&\quad\times\frac{\sqrt{2u}}{RJ_{\nu+1}(j_{\nu,n})}J_{\nu}(R^{-1}j_{\nu,n}u)\,du.
\end{aligned}
\]
Using \eqref{Volterra}, rewrite this formula as
\[
b^l_{mn}(s)=\frac{c_{NH}\sqrt{2}s^{2H-m}}{RJ_{\nu+1}(j_{\nu,n})}
\int^{\infty}_0u^{m-H}J_{\nu}(R^{-1}j_{\nu,n}u)G^{1,0}_{1,1}\left(\dfrac{u^2}{s^2}\left|
\begin{array}{l}
N/2+H\\
0
\end{array}
\right.\right)\,du.
\]
To calculate this integral,  use \eqref{22441} with $\alpha=m-H+1$, $b=R^{-1}j_{\nu,n}$, $m=p=q=1$, $r=0$, $\omega=s^{-2}$, $a_1=N/2+H$, and $b_1=0$. We obtain
\[
\begin{aligned}
b^l_{mn}(s)&=\frac{c_{NH}2^{m-H+1/2}s^{2H-m}R^{m-H}}{J_{\nu+1}(j_{\nu,n})j^{m-H+1}_{\nu,n}}\\
&\quad\times G^{1,1}_{3,1}\left(\dfrac{4R^2}{s^2j^2_{\nu,n}}\left|
\begin{array}{l}
1-\frac{m-H+1+\nu}{2},\frac{N}{2}+H,1-\frac{m-H+1-\nu}{2}\\
0
\end{array}
\right.\right).
\end{aligned}
\]

To decrease the order of the Meijer $G$-function from $4$ to $2$, put $\nu=m-1-H$ and use \eqref{8229}. We get
\[
b^l_{mn}(s)=\frac{c_{NH}2^{m-H+1/2}s^{2H-m}R^{m-H}}{J_{m-H}(j_{m-1-H,n})j^{m-H+1}_{m-1-H,n}}
G^{0,1}_{2,0}\left(\dfrac{4R^2}{s^2j^2_{m-1-H,n}}\left|
\begin{array}{l}
1-m+H,N/2+H\\
\boldsymbol{\cdot}
\end{array}
\right.\right).
\]
For further simplification, use the symmetry relation \eqref{symmetry}. We obtain
\[
b^l_{mn}(s)=\frac{c_{NH}2^{m-H+1/2}s^{2H-m}R^{m-H}}{J_{m-H}(j_{m-1-H,n})j^{m-H+1}_{m-1-H,n}}
G^{1,0}_{0,2}\left(\dfrac{s^2j^2_{m-1-H,n}}{4R^2}\left|
\begin{array}{l}
\boldsymbol{\cdot}\\
m-H,1-N/2-H\\
\end{array}
\right.\right).
\]
Then, use the argument transformation \eqref{argument} with $\alpha=m/2+1/2-N/4-H$ to get
\[
\begin{aligned}
b^l_{mn}(s)&=\frac{c_{NH}2^{H+(N-1)/2}R^{H+(N-2)/2}}{s^{(N-2)/2}J_{m-H}(j_{m-1-H,n})j^{H+N/2}_{m-1-H,n}}\\
&\quad\times G^{1,0}_{0,2}\left(\dfrac{s^2j^2_{m-1-H,n}}{4R^2}\left|
\begin{array}{l}
\boldsymbol{\cdot}\\
m/2+(N-2)/4,-m/2-(N-2)/4\\
\end{array}
\right.\right).
\end{aligned}
\]
Finally, use \eqref{bessel} to obtain
\[
\begin{aligned}
b^l_{mn}(s)&=\frac{c_{NH}2^{H+1/2}R^H}{J_{m-H}(j_{m-1-H,n})j^{H+1}_{m-1-H,n}\Gamma(N/2)}\\
&\quad\times 2^{(N-2)/2}\Gamma(N/2)\frac{J_{m+(N-2)/2}(R^{-1}j_{m-1-H,n}s)}{(R^{-1}j_{m-1-H,n}s)^{(N-2)/2}}.
\end{aligned}
\]

The last formula and \eqref{case0} can be unified as
\[
\begin{aligned}
b^l_{mn}(s)&=\frac{c_{NH}2^{H+1/2}R^H}{J_{|m-1|-H+1}(j_{|m-1|-H,n})j^{H+1}_{|m-1|-H,n}\Gamma(N/2)}\\
&\quad\times[g_m(R^{-1}j_{|m-1|-H,n}s)-\delta^m_0]
\end{aligned}
\]
with
\[
g_m(z)=2^{(N-2)/2}\Gamma(N/2)\frac{J_{m+(N-2)/2}(z)}{z^{(N-2)/2}}.
\]
Substituting the value of $c_{NH}$ from \eqref{constant}, we obtain
\[
\begin{aligned}
b^l_{mn}(s)&=\frac{2^{H+1}\sqrt{\pi^{(N-2)/2}\Gamma(N/2+H)\Gamma(H+1)\sin(\pi H)}R^H}{\Gamma(N/2)J_{|m-1|-H+1}(j_{|m-1|-H,n})j^{H+1}_{|m-1|-H,n}}\\
&\quad\times[g_m(R^{-1}j_{|m-1|-H,n}s)-\delta^m_0].
\end{aligned}
\]
This result was proved by Malyarenko (2008) for the case of $R=1$.
\end{example}

\begin{example}
For simplicity, put $R=1$. The functions
\[
\varphi_n(u)=\sqrt{2n+1}P_n(2u-1),\qquad n\geq 1,
\]
where $P_n(x)$ are Legendre polynomials, form a basis in the Hilbert space $L^2[0,1]$. The corresponding Fourier coefficients $b^l_{mn}(s)$ have the following form:
\[
\begin{aligned}
b^l_{0n}(s)&=c_{NH}\sqrt{2n+1}\int^{\infty}_0u^{H-1/2}G^{2,0}_{2,2}\left(u\left|
\begin{array}{l}
-n,n+1\\
0,0\\
\end{array}
\right.\right)G^{2,0}_{2,2}\left(\frac{u^2}{s^2}\left|
\begin{array}{l}
N/2,1\\
0,1-H\\
\end{array}
\right.\right)\,du,\\
b^l_{mn}(s)&=c_{NH}\sqrt{2n+1}s^{2H-m}\int^{\infty}_0u^{m-H-1/2}G^{2,0}_{2,2}\left(u\left|
\begin{array}{l}
-n,n+1\\
0,0\\
\end{array}
\right.\right)\\
&\quad\times G^{1,0}_{1,1}\left(\frac{u^2}{s^2}\left|
\begin{array}{l}
N/2+H\\
0\\
\end{array}
\right.\right)\,du.
\end{aligned}
\]
This follows from \eqref{coefficients} and \eqref{Gegenbauer} with $\lambda=0$.

The integrals in the last display are complicated, because the arguments of the two Meijer $G$-functions contain different powers of the independent variable $u$. However, they still can be calculated analytically, using formula 2.24.1.1 from Prudnikov et al (1990). The answer is
\[
\begin{aligned}
b^l_{0n}(s)&=\frac{c_{NH}\sqrt{2n+1}}{2}\\
&\quad\times G^{2,4}_{6,6}\left(s^{-2}\left|
\begin{array}{l}
\frac{1-2H}{4},\frac{3-2H}{4},\frac{1-2H}{4},\frac{3-2H}{4},\frac{N}{2},1\\
0,1-H,\frac{1-2H+2n}{4},\frac{3-2H+2n}{4},\frac{-1-2H-2n}{4},\frac{1-2H-2n}{4}\\
\end{array}
\right.\right),\\
b^l_{mn}(s)&=\frac{c_{NH}\sqrt{2n+1}s^{2H-m}}{2}\\
&\quad\times G^{1,4}_{5,5}\left(s^{-2}\left|
\begin{array}{l}
\frac{1-2m+2H}{4},\frac{3-2m+2H}{4},\frac{1-2m+2H}{4},\frac{3-2m+2H}{4},\frac{N}{2}+H\\
0,\frac{1-2H+2n-2m}{4},\frac{3-2H+2n+2m}{4},\frac{-1-2H-2n-2m}{4},\frac{1-2H-2n-2m}{4}\\
\end{array}
\right.\right).
\end{aligned}
\]
The details are left to the reader.
\end{example}

\section{Proofs of technical lemmas}\label{technical}

\begin{proof}[Proof of Lemma~\ref{l1}]
It follows from \eqref{auto} that the autocorrelation function of the multiparameter fractional Brownian motion can be written as
\[
R(s,t,u)=\frac{1}{2}(s^{2H}+t^{2H}-(s^2-2stu+t^2)^H).
\]
Therefore, the autocorrelation function of the stochastic process $X^l_m(t)$ has the form
\[
R_m(s,t)=I_1-\lim_{\alpha\to(N-1)/2}I_2(\alpha),
\]
where
\[
\begin{aligned}
I_1&=\frac{c}{2}(s^{2H}+t^{2H})\int^1_{-1}C^{(N-2)/2}_m(u)(1-u^2)^{(N-3)/2}\,du,\\
I_2(\alpha)&=\frac{c}{2}(2st)^H\int^1_{-1}(u+1)^{\alpha-1}(1-u)^{(N-3)/2}
\left(\frac{s^2+t^2}{2st}-u\right)^HC^{(N-2)/2}_m(u)\,du.
\end{aligned}
\]

To calculate $I_1$, we use formula 2.21.2.17 from Prudnikov et al (1988):
\[
\int^a_{-a}(a^2-x^2)^{\lambda-1/2}C^{\lambda}_n(x/a)\,dx=\delta^0_m\sqrt{\pi}a^{2\lambda}
\frac{\Gamma(\lambda+1/2)}{\Gamma(\lambda+1)}
\]
with $a=1$, $\lambda=(N-2)/2$, and $n=m$. Here, $\delta^0_m$ is the Kronecker's delta. After simplification, we obtain
\[
I_1=\frac{2^{N-3}\pi^{(N-1)/2}\Gamma((N-2)/2)\Gamma((N-1)/2)}{\Gamma(N-2)\Gamma(N/2)}
(s^{2H}+t^{2H})\delta^0_m.
\]
This expression can be further simplified using the doubling formula\footnote{\href{http://functions.wolfram.com/06.05.16.0006.01}{\texttt{http://functions.wolfram.com/06.02.16.0006.01}}}
\[
\Gamma(2z)=\frac{2^{2z-1}}{\sqrt{\pi}}\Gamma(z)\Gamma(z+1/2)
\]
with $z=(N-2)/2$. We get
\[
I_1=\frac{\pi^{N/2}}{\Gamma(N/2)}(s^{2H}+t^{2H})\delta^0_m.
\]

To calculate $I_2(\alpha)$, we use formula 2.21.4.15 from Prudnikov et al (1988):
\[
\begin{aligned}
&\int^a_{-a}(x+a)^{\alpha-1}(a-x)^{\lambda-1/2}
(z-x)^{-\vartheta}C^{\lambda}_n(x/a)\,dx=\\
&\quad\frac{(-1)^n}{n!}(1/2+\lambda-\alpha)_n(2\lambda)_n\frac{\Gamma(\alpha)\Gamma(\lambda+1/2)}
{\Gamma(\alpha+\lambda+n+1/2)}(2a)^{\alpha+\lambda-1/2}(z+a)^{-\vartheta}\\
&\qquad\times{}_3F_2(\alpha,\theta,1/2+\alpha-\lambda;1/2+\alpha-\lambda-n,1/2+\alpha+\lambda+n;
2a/(a+z))
\end{aligned}
\]
with $a=1$, $\lambda=(N-2)/2$, $z=(s^2+t^2)/(2st)$, $\vartheta=-H$, and $n=m$. After simplification, we obtain
\[
\begin{aligned}
I_2(\alpha)&=(-1)^m2^{\alpha+3(N-3)/2}\pi^{(N-2)/2}\\
&\quad\times\frac{\Gamma((N-2)/2)\Gamma((N-1)/2-\alpha+m)\Gamma(\alpha)\Gamma((N-1)/2)}
{\Gamma((N-1)/2-\alpha)\Gamma(N-2)\Gamma(\alpha+(N-1)/2+m)}(s+t)^{2H}\\
&\quad\times{}_3F_2\left(\alpha,-H,\alpha-\frac{N-3}{2};\alpha-\frac{N-3}{2}-m,\alpha+\frac{N-1}{2}+m;
\frac{4st}{(s+t)^2}\right).
\end{aligned}
\]
Using the doubling formula with $z=(N-2)/2$, we get
\begin{equation}\label{eq3}
\begin{aligned}
I_2(\alpha)&=(-1)^m2^{\alpha+(N-3)/2}\pi^{(N-1)/2}\frac{\Gamma((N-1)/2-\alpha+m)\Gamma(\alpha)(s+t)^{2H}}
{\Gamma((N-1)/2-\alpha)\Gamma(\alpha+(N-1)/2+m)}\\
&\quad\times{}_3F_2\left(\alpha,-H,\alpha-\frac{N-3}{2};\alpha-\frac{N-3}{2}-m,\alpha+\frac{N-1}{2}+m;
\frac{4st}{(s+t)^2}\right).
\end{aligned}
\end{equation}

In the case of $m=0$, \eqref{eq3} simplifies as follows.
\[
\begin{aligned}
I_2(\alpha)&=2^{\alpha+(N-3)/2}\pi^{(N-1)/2}\frac{\Gamma(\alpha)(s+t)^{2H}}
{\Gamma(\alpha+(N-1)/2)}\\
&\quad\times{}_3F_2\left(\alpha,-H,\alpha-\frac{N-3}{2};\alpha-\frac{N-3}{2},\alpha+\frac{N-1}{2};
\frac{4st}{(s+t)^2}\right).
\end{aligned}
\]
According to paragraph 7.2.3.2 from Prudnikov et al (1990), the value of the generalised hypergeometric function ${}_pF_q(a_1,\dots,a_p;b_1,\dots,b_q;z)$ is independent on the order of \emph{upper parameters} $a_1$, \dots, $a_p$ and \emph{lower parameters} $b_1$, \dots, $b_q$. Moreover, formula 7.2.3.7 from Prudnikov et al (1990) states that
\begin{equation}\label{decrease}
\begin{aligned}
&{}_pF_q(a_1,\dots,a_{p-r},c_1,\dots,c_r;b_1,\dots,b_{q-r},c_1,\dots,c_r;z)\\
&\quad={}_{p-r}F_{q-r}(a_1,\dots,a_{p-r};b_1,\dots,b_{q-r};z).
\end{aligned}
\end{equation}
Using these properties, we get
\[
I_2(\alpha)=2^{\alpha+(N-3)/2}\pi^{(N-1)/2}\frac{\Gamma(\alpha)(s+t)^{2H}}
{\Gamma(\alpha+(N-1)/2)}{}_2F_1\left(\alpha,-H;\alpha+\frac{N-1}{2};\frac{4st}{(s+t)^2}\right).
\]
In particular,
\[
\begin{aligned}
\lim_{\alpha\to(N-1)/2}I_2(\alpha)&=2^{N-2}\pi^{(N-1)/2}\frac{\Gamma((N-1)/2)(s+t)^{2H}}
{\Gamma(N-1)}\\
&\quad\times{}_2F_1\left(\frac{N-1}{2},-H;N-1;\frac{4st}{(s+t)^2}\right).
\end{aligned}
\]
The application of the doubling formula with $z=(N-1)/2$ yields
\[
\lim_{\alpha\to(N-1)/2}I_2(\alpha)=\frac{\pi^{N/2}}{\Gamma(N/2)}(s+t)^{2H}
{}_2F_1\left(\frac{N-1}{2},-H;N-1;\frac{4st}{(s+t)^2}\right).
\]
The argument simplification formula\footnote{\href{http://functions.wolfram.com/07.23.16.0005.01}{\texttt{http://functions.wolfram.com/07.23.16.0005.01}}} states:
\[
{}_2F_1(a,b;2b;4z/(z+1)^2)=(z+1)^{2a}{}_2F_1(a,a-b+1/2;b+1/2;z^2).
\]
Use this formula with $a=-H$, $b=(N-1)/2$, and $z=s/t$. We obtain
\[
\lim_{\alpha\to(N-1)/2}I_2(\alpha)=\frac{\pi^{N/2}}{\Gamma(N/2)}t^{2H}{}_2F_1(-H,1-H-N/2;N/2;s^2/t^2).
\]
and, finally
\begin{equation}\label{r0}
R_0(s,t)=\frac{\pi^{N/2}}{\Gamma(N/2)}[s^{2H}+t^{2H}(1-{}_2F_1(-H,1-H-N/2;N/2;s^2/t^2))].
\end{equation}

It remains to prove that the first case in \eqref{xlmauto} simplifies to \eqref{r0}. To do this, use the representation of the Meijer $G$-function through hypergeometric functions\footnote{\href{http://functions.wolfram.com/07.34.26.0004.01}{\texttt{http://functions.wolfram.com/07.34.26.0004.01}}} \[
\begin{aligned}
&G^{m,n}_{p,q}\left(z\left|
\begin{array}{l}
a_1,\dots,a_p\\
b_1,\dots,b_q
\end{array}
\right.\right)=\sum^m_{k=1}\frac{\prod_{j\in\{1,2,\dots,m\}\setminus\{k\}}\Gamma(b_j-b_k)
\prod_{j=1}^m\Gamma(1-a_j+b_k)}{\prod_{j=n+1}^p\Gamma(a_j-b_k)\prod_{j=m+1}^q\Gamma(1-b_j+b_k)}z^{b_k}\\
&\quad\times{}_pF_{q-1}\left(
\begin{array}{l}
1-a_1+b_k,\dots,1-a_p+b_k;\\
1+a_1-a_k,\dots,1+a_{k-1}-a_k,1+a_{k+1}-a_k,\dots,
1+a_q-a_k;(-1)^{p-m-n}z
\end{array}
\right)
\end{aligned}
\]
with $m=n=2$, $p=q=4$, $z=s^2/t^2$, $a_1=1-H$, $a_2=b_1=0$, $a_3=N/2$, $a_4=1$, $b_2=1-H$, $b_3=1-H-N/2$, and $b_4=-H$. We obtain
\[
\begin{aligned}
I_1&=\frac{c_{NH}^2}{2}s^{2H}\left[\frac{\Gamma(1-H)\Gamma(H)}{\Gamma(N/2)\Gamma(N/2+H)\Gamma(1+H)}
{}_4F_3\left(
\begin{array}{l}
H,1,1-N/2,0;\\
H,H+N/2,1+H;s^2/t^2
\end{array}
\right)\right.\\
&\quad\left.+\frac{\Gamma(H-1)\Gamma(2-H)}{\Gamma(N/2+H-1)\Gamma(H)\Gamma(N/2+1)}
\frac{s^{2(1-H)}}{t^{2(1-H)}}{}_4F_3\left(
\begin{array}{l}
1,2-H,2-H-N/2,1-H;\\
2-H,N/2+1,2;s^2/t^2
\end{array}
\right)\right].
\end{aligned}
\]

The first term is simplified, using the following formula\footnote{\href{http://functions.wolfram.com/07.31.03.0012.01}{\texttt{http://functions.wolfram.com/07.31.03.0012.01}}}
\[
{}_pF_q(0,a_2,\dots,a_p;b_1,\dots,b_q;z)=1,
\]
while the second term is simplified by \eqref{decrease}. We get
\[
\begin{aligned}
I_1&=\frac{c_{NH}^2}{2}s^{2H}\left[\frac{\Gamma(1-H)\Gamma(H)}{\Gamma(N/2)\Gamma(N/2+H)\Gamma(1+H)}\right.\\
&\quad\left.+\frac{\Gamma(H-1)\Gamma(2-H)}{\Gamma(N/2+H-1)\Gamma(H)\Gamma(N/2+1)}
\frac{s^{2(1-H)}}{t^{2(1-H)}}{}_3F_2\left(
\begin{array}{l}
1,2-H-N/2,1-H;\\
N/2+1,2;s^2/t^2
\end{array}
\right)\right].
\end{aligned}
\]
Using the formula\footnote{\href{http://functions.wolfram.com/07.27.03.0120.01}{\texttt{http://functions.wolfram.com/07.27.03.0120.01}}}
\[
{}_3F_2(1,b,c;2,e;z)=\frac{e-1}{(b-1)(c-1)z}[{}_2F_1(b-1,c-1;e-1;z)-1].
\]
with $b=2-H-N/2$, $c=1-H$, $e=N/2+1$, and $z=s^2/t^2$ yields
\[
\begin{aligned}
I_1&=\frac{c_{NH}^2}{2}s^{2H}\left[\frac{\Gamma(1-H)\Gamma(H)}{\Gamma(N/2)\Gamma(N/2+H)\Gamma(1+H)}\right.\\
&\quad+\frac{\Gamma(H-1)\Gamma(2-H)(N/2)}{\Gamma(N/2+H-1)\Gamma(H)\Gamma(N/2+1)(1-H-N/2)(-H)}
\frac{t^{2H}}{s^{2H}}\\
&\quad\times\left.[{}_2F_1(1-H-N/2,-H;N/2;s^2/t^2)-1]\right].
\end{aligned}
\]
To simplify the second line, use \eqref{recurrent} in the following way:
\[
\begin{aligned}
\Gamma(H-1)&=\Gamma(H)/((H-1),\\
\Gamma(2-H)&=(1-H)\Gamma(1-H),\\ \Gamma(N/2+H-1)(1-H-N/2)&=-\Gamma(N/2+H),\\ \Gamma(H)(-H)&=-\Gamma(H+1),\\
\Gamma(N/2+1)&=(N/2)\Gamma(N/2).
\end{aligned}
\]
After simplification, we obtain
\[
I_1=\frac{c_{NH}^2\Gamma(1-H)\Gamma(H)}{2\Gamma(N/2)\Gamma(N/2+H)\Gamma(H+1)}
[s^{2H}+t^{2H}(1-{}_2F_1(-H,1-H-N/2;N/2;s^2/t^2))],
\]
or, by \eqref{additional} with $z=H$,
\[
\begin{aligned}
I_1&=\frac{c_{NH}^2\pi}{2\Gamma(N/2)\Gamma(N/2+H)\Gamma(H+1)\sin(\pi H)}\\
&\quad\times[s^{2H}+t^{2H}(1-{}_2F_1(-H,1-H-N/2;N/2;s^2/t^2))].
\end{aligned}
\]
Substituting \eqref{constant} to the last display, we get
\[
I_1=\frac{\pi^{N/2}}{\Gamma(N/2)}[s^{2H}+t^{2H}(1-{}_2F_1(-H,1-H-N/2;N/2;s^2/t^2))].
\]
This completes the proof of the case of $m=0$.

In the case of $m\geq 1$, \eqref{eq3} can be rewritten as
\[
\begin{aligned}
\lim_{\alpha\to(N-1)/2}I_2(\alpha)&=(-1)^m2^{N-2}\pi^{(N-1)/2}\frac{(m-1)!\Gamma((N-1)/2)(s+t)^{2H}}
{\Gamma(N-1+m)}\\
&\quad\times\lim_{\alpha\to(N-1)/2}\frac{{}_3F_2\left(\alpha,-H,\alpha-\frac{N-3}{2};
\alpha-\frac{N-3}{2}-m,N-1+m;\frac{4st}{(s+t)^2}\right)}{\Gamma(\alpha-(N-3)/2-m)}\\
&\quad\times\lim_{\alpha\to(N-1)/2}\frac{\Gamma(\alpha-(N-3)/2-m)}{\Gamma((N-1)/2-\alpha)}.
\end{aligned}
\]
To calculate the first limit, we use the following formula\footnote{\href{http://functions.wolfram.com/07.31.25.0003.01}{\texttt{http://functions.wolfram.com/07.31.25.0003.01}}}
\[
\begin{aligned}
&\lim_{b_1\to-n}\frac{{}_pF_q(a_1,\dots,a_p;b_1,\dots,b_q;z)}{\Gamma(b_1)}=
\frac{z^{n+1}}{(n+1)!}\frac{\prod^p_{j=1}(a_j)_{n+1}}{\prod^q_{j=2}(b_j)_{n+1}}\\
&\quad\times{}_pF_q(a_1+n+1,\dots,a_p+n+1;b_2+n+1,\dots,b_q+n+1,n+2;z),\quad n\in\mathbb{Z}_+
\end{aligned}
\]
with $n=m-1$, $p=3$, $q=2$, $a_1=(N-1)/2$, $a_2=-H$, $a_3=1$, $b_1=1-m$, $b_2=N-1+m$, and $z=4st/(s+t)^2$. We get
\[
\begin{aligned}
&\lim_{\alpha\to(N-1)/2}\frac{{}_3F_2\left(\alpha,-H,\alpha-\frac{N-3}{2};
\alpha-\frac{N-3}{2}-m,N-1+m;\frac{4st}{(s+t)^2}\right)}{\Gamma(\alpha-(N-3)/2-m)}=\\
&\quad\frac{2^{2m}(st)^m\Gamma((N-1)/2+m)\Gamma(m-H)\Gamma(N-1+m)}
{(s+t)^{2m}\Gamma((N-1)/2)\Gamma(-H)\Gamma(N-1+2m)}\\
&\qquad\times{}_3F_2((N-1)/2+m,m-H,m+1;N-1+2m,m+1;4st/(s+t)^2).
\end{aligned}
\]
Using the doubling formula with $z=(N-1)/2+m$ and \eqref{decrease}, we get
\[
\begin{aligned}
&\lim_{\alpha\to(N-1)/2}\frac{{}_3F_2\left(\alpha,-H,\alpha-\frac{N-3}{2};
\alpha-\frac{N-3}{2}-m,N-1+m;\frac{4st}{(s+t)^2}\right)}{\Gamma(\alpha-(N-3)/2-m)}=\\
&\quad\frac{(st)^m\Gamma(m-H)\Gamma(N-1+m)\sqrt{\pi}}
{(s+t)^{2m}\Gamma((N-1)/2)\Gamma(-H)2^{N-2}\Gamma(N/2+m)}\\
&\qquad\times{}_2F_1((N-1)/2+m,m-H;N-1+2m;4st/(s+t)^2).
\end{aligned}
\]

Rewrite the second limit as
\[
\begin{aligned}
\lim_{\alpha\to(N-1)/2}\frac{\Gamma(\alpha-(N-3)/2-m)}{\Gamma((N-1)/2-\alpha)}&=
\lim_{\alpha\to(N-1)/2}\frac{\Gamma(\alpha-(N-3)/2-m)}{\Gamma(\alpha-(N-1)/2}\\
&\quad\times\lim_{\alpha\to(N-1)/2}\frac{\Gamma(\alpha-(N-1)/2)}{\Gamma((N-1)/2-\alpha)}.
\end{aligned}
\]
For the first part, use the formula\footnote{\href{http://functions.wolfram.com/06.05.16.0022.01}{\texttt{http://functions.wolfram.com/06.05.16.0022.01}}}
\[
\frac{\Gamma(z-n)}{\Gamma(z)}=\prod^n_{k=1}\frac{1}{z-k}
\]
with $z=\alpha-(N-1)/2$ and $n=m-1$. For the second part, use \eqref{recurrent} with $z=\alpha-(N-3)/2$. We get
\[
\lim_{\alpha\to(N-1)/2}\frac{\Gamma(\alpha-(N-3)/2-m)}{\Gamma((N-1)/2-\alpha)}=
\frac{(-1)^{m-1}}{(m-1)!}(-1)=\frac{(-1)^m}{(m-1)!}.
\]

Combining everything together, we obtain
\[
\begin{aligned}
R_m(s,t)&=-\frac{\pi^{N/2}(s+t)^{2(H-m)}(st)^m\Gamma(m-H)}{\Gamma(-H)\Gamma(N/2+m)}\\
&\quad\times{}_2F_1((N-1)/2+m,m-H;N-1+2m;4st/(s+t)^2).
\end{aligned}
\]
Use \eqref{recurrent} with $z=1-H$:
\[
\begin{aligned}
R_m(s,t)&=\frac{\pi^{N/2}(s+t)^{2(H-m)}(st)^m\Gamma(m-H)H}{\Gamma(1-H)\Gamma(N/2+m)}\\
&\quad\times{}_2F_1((N-1)/2+m,m-H;N-1+2m;4st/(s+t)^2).
\end{aligned}
\]
By the additional formula for the gamma function\footnote{\href{http://functions.wolfram.com/06.05.16.0010.01}{\texttt{http://functions.wolfram.com/06.05.16.0010.01}}}
\begin{equation}\label{additional}
\Gamma(z)\Gamma(1-z)=\frac{\pi}{\sin(\pi z)}
\end{equation}
with $z=H$, and \eqref{recurrent} with $z=H+1$ we get
\[
\begin{aligned}
R_m(s,t)&=\frac{\pi^{(N-2)/2}(s+t)^{2(H-m)}(st)^m\Gamma(m-H)\Gamma(H+1)\sin(\pi H)}{\Gamma(N/2+m)}\\
&\quad\times{}_2F_1((N-1)/2+m,m-H;N-1+2m;4st/(s+t)^2),
\end{aligned}
\]
and then, using the following formula\footnote{\href{http://functions.wolfram.com/07.23.26.0031.01}{\texttt{http://functions.wolfram.com/07.03.26.0031.01}}} \[
\begin{aligned}
(\sqrt{z}+1)^{-2a}{}_2F_1(a,b;2b;4\sqrt{z}/(\sqrt{z}+1)^2)&=\frac{\Gamma(b+1/2)\Gamma(b-a+1/2)}{\Gamma(a)}\\
&\quad\times G^{1,1}_{2,2}\left(z\left|
\begin{array}{l}
1-a,b-a+1/2\\
0,1/2-b
\end{array}
\right.\right)
\end{aligned}
\]
with $z=s^2/t^2$, $a=m-H$, and $b=(N-1)/2+m$, we finally obtain
\[
R_m(s,t)=\frac{c^2_{NH}}{2}s^mt^{2H-m}G^{1,1}_{2,2}\left(\dfrac{s^2}{t^2}\left|
\begin{array}{l}
H+1-m,N/2+H\\
0,1-N/2-m
\end{array}
\right.\right),\quad m\geq 1.
\]
This completes the proof.
\end{proof}

\begin{proof}[Proof of Lemma~\ref{l2}]
\eqref{Volterra} for the case of $m=0$ is obvious from the formula\footnote{\href{http://functions.wolfram.com/07.34.03.0645.01}{\texttt{http://functions.wolfram.com/07.34.03.0645.01}}}
\[
\begin{aligned}
G^{2,0}_{2,2}\left(x\left|
\begin{array}{l}
a_1,a_2\\
b_1,b_2
\end{array}
\right.\right)&=\frac{x^{b_1}(1-x)^{a_1+a_2-b_1-b_2-1}_+}{\Gamma(a_1+a_2-b_1-b_2)}\\
&\quad\times{}_2F_1(a_2-b_2,a_1-b_2;a_1+a_2-b_1-b_2;1-x)
\end{aligned}
\]
with $x=u^2/s^2$, $a_1=N/2$, $a_2=1$, $b_1=0$, and $b_2=1-H$, where we use a shortcut $(1-x)_+$ for $\max\{1-x,0\}$. Similarly, \eqref{Volterra} for the case of $m\geq 1$ is obvious from the formula\footnote{\href{http://functions.wolfram.com/07.34.03.0247.01}{\texttt{http://functions.wolfram.com/07.34.03.0247.01}}} \[
G^{1,0}_{1,1}\left(x\left|
\begin{array}{l}
a\\
b
\end{array}
\right.\right)=\frac{x^b(1-x)^{a-b-1}_+}{\Gamma(a-b)}
\]
with $x=u^2/s^2$, $a=N/2+H$, and $b=0$.

It remains to calculate two integrals. The first one is as follows:
\[
I_1=c_{NH}^2\int^{\min\{s,t\}}_0u^{2H-1}G^{2,0}_{2,2}\left(\dfrac{u^2}{s^2}\left|
\begin{array}{l}
N/2,1\\
0,1-H
\end{array}
\right.\right)G^{2,0}_{2,2}\left(\dfrac{u^2}{t^2}\left|
\begin{array}{l}
N/2,1\\
0,1-H
\end{array}
\right.\right)\,du.
\]
Taking into account \eqref{Volterra}, we can substitute $\infty$ for the upper limit of integration. After change of variable $u=\sqrt{x}$, we obtain
\[
I_1=\frac{c_{NH}^2}{2}\int^{\infty}_0x^{H-1}G^{2,0}_{2,2}\left(\dfrac{x}{s^2}\left|
\begin{array}{l}
N/2,1\\
0,1-H
\end{array}
\right.\right)G^{2,0}_{2,2}\left(\dfrac{x}{t^2}\left|
\begin{array}{l}
N/2,1\\
0,1-H
\end{array}
\right.\right)\,dx.
\]

To calculate this integral, use \eqref{classic} with $\alpha=H$, $s=m=u=v=p=q=2$, $t=v=0$, $w=1/s^2$, $z=1/t^2$, $c_1=a_1=N/2$, $c_2=a_2=1$, $d_1=b_1=0$, and $d_2=b_2=1-H$. We get
\[
I_1=\frac{c_{NH}^2}{2}s^{2H}G^{2,2}_{4,4}\left(\dfrac{s^2}{t^2}\left|
\begin{array}{l}
1-H,0,N/2,1\\
0,1-H,1-H-N/2,-H
\end{array}
\right.\right).
\]
This completes the calculation of the first integral.

The second integral is as follows.
\[
I_2=c_{NH}^2(st)^{2H-m}\int^{\min\{s,t\}}_0 u^{2m-2H-1}G^{1,0}_{1,1}\left(\dfrac{u^2}{s^2}\left|
\begin{array}{l}
N/2+H\\
0
\end{array}
\right.\right)G^{1,0}_{1,1}\left(\dfrac{u^2}{t^2}\left|
\begin{array}{l}
N/2+H\\
0
\end{array}
\right.\right)\,du.
\]
Taking into account \eqref{Volterra}, we can substitute $\infty$ for the upper limit of integration. After change of variable $u=\sqrt{x}$, we obtain
\[
I_2=\frac{c_{NH}^2}{2}(st)^{2H-m}\int^{\infty}_0 x^{m-H-1}G^{1,0}_{1,1}\left(\dfrac{x}{s^2}\left|
\begin{array}{l}
N/2+H\\
0
\end{array}
\right.\right)G^{1,0}_{1,1}\left(\dfrac{x}{t^2}\left|
\begin{array}{l}
N/2+H\\
0
\end{array}
\right.\right)\,dx.
\]
To calculate this integral, use \eqref{classic} with $\alpha=m-H$, $s=m=u=v=p=q=1$, $t=v=0$, $w=1/s^2$, $z=1/t^2$, $c_1=a_1=N/2+H$, and $d_1=b_1=0$. We get
\[
I_2=\frac{c_{NH}^2}{2}s^mt^{2H-m}G^{1,1}_{2,2}\left(\dfrac{s^2}{t^2}\left|
\begin{array}{l}
H+1-m,N/2+H\\
0,1-N/2-H
\end{array}
\right.\right).
\]
This completes the proof.
\end{proof}


\begin{thebibliography}{00}

\bibitem{Dec} Decreusefond, L. (2005). Stochastic integration with respect to Volterra processes, \emph{Ann. I. H.~Poincar\'{e}}, \textbf{41}, 123--149.

\bibitem{Kol} Kolmogorov, A. (1940). Wienerische Spiralen und einige andere interessante Kurven im Hilbertschen Raum.  \emph{C.R. (Doklady) Acad. Sci. USSR (N.S.)}, \textbf{26}, 115--118.

\bibitem{Mal} Malyarenko, A. (2008). An Optimal Series Expansion of the Multiparameter Fractional Brownian Motion. \emph{J. Theor. Probab.}, \textbf{21}, 459--475.

\bibitem{Man} Mandelbrot, B., and J. van Ness (1968). Fractional Brownian motions, fractional noises and applications, \emph{SIAM Rev.}, \textbf{10}, 422--437.

\bibitem{Pru2} Prudnikov, A.P., Brychkov, Yu.A., and O.I. Marichev (1988). \emph{Integrals and series. Vol. 2. Special functions}, Second edition, Gordon \& Breach Science Publishers, New York.

\bibitem{Pru3} Prudnikov, A.P., Brychkov, Yu.A., and O.I. Marichev (1990). \emph{Integrals and series. Vol. 3. More special functions}, Gordon \& Breach Science Publishers, New York.

\bibitem{Wat} Watson, G.N. (1995). \emph{A treatise on the
    theory of Bessel functions}, Cambridge University Press, Cambridge.

\bibitem{Yad} Yadrenko, M.I. (1983). \emph{Spectral theory of random fields}, Optimization Software, New York.

\end{thebibliography}
\end{document}